\documentclass{article}
\usepackage{amsmath,amsfonts,amssymb,amscd,amsthm}
\usepackage{psfrag,extarrows,multirow}
\usepackage{color,graphics}
\usepackage{mathdots}
\usepackage[english]{babel}

\usepackage{tikz,tikz-3dplot}
\usetikzlibrary{matrix,arrows,automata,intersections}

\newtheorem{thm}{Theorem}
\newtheorem{de}{Definition}

\newtheorem{coro}{Corollary}
\newtheorem{lem}{Lemma}
\newcommand{\be}{\begin{equation}}
\newcommand{\ee}{\end{equation}}

\begin{document}
\title{K\"ahler Finsler Metrics and Conformal Deformations\footnote{Supported by the National Natural Science Foundation of China (no. 11871126, 11471246,  11101307).}}
\author{Bin Chen\ \ Yibing Shen\ \  Lili Zhao}
\maketitle

\begin{abstract}
The conformal properties of complex Finsler metrics are studied.  We give a characterization  of a compact complex Finsler manifold to be globally conformal K\"ahler. The critical points of the total holomorphic curvature and total Ricci curvature in the volume preserved conformal classes are studied. The stability of critical K\"ahler Finsler metrics is obtained. A Yamabe type problem for mean Ricci curvature is considered.


\end{abstract}

\section{Introduction}\setcounter{equation}{0}\setcounter{thm}{0}\setcounter{de}{0}\setcounter{prop}{0}\setcounter{coro}{0}\setcounter{lem}{0}

Searching for the notion of the ``best" metric on a manifold is a central problem in geometry. In Riemannian realm, the canonical ones are Yamabe metrics, Einstein metrics and etc. In complex geometry, one is led to extremal metrics, K\"ahler Einstein metrics and etc. During the past decades, there is a bundle of results on the ``best" Finsler metrics, such as Einstein Finsler metrics, Yamabe Finsler metrics and etc. (cf. \cite{AP,BR,CZ} and references therein). Complex Finsler metrics are natural generalization of Hermitian metrics. Since the most often used intrinsic (depending only on the complex structure) metrics are generally
Finsler ones (such as Kobayashi metric and Carath\'{e}odory metric), it is one hot issue to develop the theory of complex Finsler geometry. In this paper, we will study some canonical  complex Finsler metrics in a conformal class. The manifolds considered in this paper are of the complex dimension $n\geq 2$.

 The concept of K\"ahler Finsler metrics is introduced by M. Abate and G. Patrizio  in \cite{AP}. The global properties of K\"ahler Finsler spaces are well studied. The Hodge decomposition theorem  is proved  by C. Zhong and T. Zhong \cite{ZZ}. Later, J. Han and the second author study the existence of harmonic maps \cite{HS}. Recently,  the comparison theorems are obtained by S. Yin and X. Zhang \cite{YZ}.

 The first goal of this paper is to study the existence of K\"ahler Finsler metrics in a conformal class.
 Let $M$ be an $n$-dimensional compact complex space with a complex Finsler metric $G$, whose volume preserved conformal class is denoted by $[G]$. It is natural to ask whether there exists  a K\"ahler Finsler metric in $[G]$. The uniqueness is easy to obtain.

\begin{thm} In the volume preserved conformal class $[G]$, there exists at most one K\"ahler Finsler metric.
\end{thm}

In order to get the existence of K\"ahler Finsler metrics in $[G]$,  we should work on \textit{K\"ahler Finsler manifolds}. A manifold $M$ is called a K\"ahler Finsler manifold if it admits a  K\"ahler Finsler metric.
 \begin{thm} Let $M$ be a compact  K\"ahler Finsler manifold, and $G$ be an arbitrary complex Finsler metric (not necessarily K\"ahlerian) on $M$. Then, there exists a K\"ahler Finsler metric in $[G]$  if and only if the horizontal torsion of $G$ is reducible and the real part of its mean horizontal torsion is closed.
\end{thm}
The exact meaning of reducibility of the horizontal torsion can be found in Theorem 4.3.

The second goal of this paper is to understand the curvature behavior of a K\"ahler Finsler metric in its conformal class. Applying the integration along the fibre of the projectivized tangent bundle over $M$, we introduce
 the \textit{mean holomorphic curvature} $\kappa=\kappa(z)$  (see (5.18)) and  the \textit{mean holomorphic Ricci curvature} $\rho=\rho(z)$ (see (6.8)). By considering the following two total curvature functionals
\be \mathcal{K}(G)=\int_M\kappa\,d\mu_M,\ \ \mathcal{R}(G)=\int_M\rho\,d\mu_M,\ee
we obtain the following result.

\begin{thm} Let $G$ be a  K\"ahler Finsler metric on a compact complex manifold.

\noindent $(i)$  $G$ is a critical point of $\mathcal{K}$ in $[G]$ if and only if $\kappa =const$.  Moreover, $G$  is stable if and only if
  $\kappa\leq\lambda_1^h$.

\noindent $(ii)$  $G$ is a critical point of $\mathcal{R}$ in $[G]$ if and only if $\rho =const$.
Moreover, $G$  is stable if and only if  $\rho\leq\lambda_1^g$.
\end{thm}
Here $\lambda_1^h$ and $\lambda_1^g$ are the first eigenvalues of the Hermitian Laplacian of the metric measure spaces $(M,h,d\mu_M)$ and $(M,g,d\mu_M)$ respectively, where the induced metrics $h$ and $g$ are given by (5.29) and (6.14).
We shall remark that  the total holomorphic curvature was firstly considered by J. Bland and M. Kalka and the variation formula was obtained in \cite{BK}.

A K\"ahler Finsler metric is said to be \textit{Einstein} if its holomorphic Ricci curvature is constant. One can immediately get the following corollary.

\begin{coro}A K\"ahler Einstein Finsler metric  with non-positive holomorphic Ricci curvature is a stable critical point of $\mathcal{R}$ in its volume preserved conformal class.
\end{coro}

The last goal of this paper is to consider a Yamabe type problem. For a complex Finsler metric which is not necessarily K\"ahlerian,  the \textit{$\vartheta$-mean holomorphic Ricci curvature} $\rho_\vartheta$ is introduced (see (6.9)) . We then study the existence of conformal metrics with constant $\rho_\vartheta$. In the real Finsler geometry,  a similar problem is considered in \cite{CZ} for ``C-convex" metrics. It is interesting that the C-convexity  is not needed in the complex realm. Precisely, by introducing the conformal invariants $Y(G)$  and $C(G)$ (see (7.4) and (7.11) respectively), we prove the following existence theorem.

\begin{thm} Let $(M,G)$ be a compact complex Finsler manifold with complex dimension $n$. It always holds $Y(G)\cdot C(G)\leq\frac{\sigma_{2n}}{2n-2}$  where $\sigma_{2n}$ is the best Sobolev constant. If $Y(G)\cdot C(G)<\frac{\sigma_{2n}}{2n-2}$, then there exists a metric  with constant $\rho_\vartheta$ in the conformal class $[G]$.
\end{thm}

The contents of this paper are arranged as follows. In \S2, we give a brief overview
of complex Finsler metrics and the K\"ahler condition. In \S3, we introduce the integration along the fibre of the projectivized tangent bundle. In \S4, the notions of locally conformal K\"ahler and globally conformal K\"ahler are given, and Theorem 1.1 and 1.2 are proved. In \S5, we consider the functional $\mathcal{K}$ and obtain the first part of Theorem 1.3. In \S6, the functional $\mathcal{R}$ is studied and the second part of Theorem 1.3 is obtained. In the last section, the Yamabe type problem is considered and Theorem 1.4 is verified.

\section{Complex Finsler metrics}\setcounter{equation}{0}\setcounter{thm}{0}\setcounter{de}{0}\setcounter{prop}{0}\setcounter{coro}{0}\setcounter{lem}{0}

Let $M$ be a complex manifold with $\dim_\mathbb{C}M=n$, and $T'M$ be the holomorphic tangent bundle. The points of $T'M$ will be denoted by $(z,v)$ where $v=v^i\partial/\partial z^i\in T_z'M$, and thus $(z^i;v^i)$ forms a local holomorphic coordinate system of $T'M$.  Let us denote the slit holomorphic tangent bundle $T'M\setminus\{\mathbf{0}\}$ by $\tilde M$. A \textit{complex Finsler metric} on $M$ is a continuous function $G:T'M\to [0,+\infty)$ satisfies

(I)  $G(z,v)\geq 0$, where the equality holds if and only if
$v=0$;

(II) $G(z,v)\in
C^{\infty}(\tilde{M})$;

(III) $G(z,\lambda v)=\lambda\bar\lambda G(z,v)$ for $\lambda \in
\mathbb{C}^*=\mathbb{C}\setminus\{0\}$;

(IV) the Levi matrix $(G_{i\bar j})_{n\times n}:=\left(\frac{\partial^2G}{\partial v^i\partial\bar v^j}\right)_{n\times n}$ is positively definite on $\tilde M$.

\noindent
\noindent The last condition is called the  {\it strongly pseudo-convexity} of $G$.   The pair $(M,G)$ is called a \textit{complex Finsler manifold}.
Throughout this paper, all the manifolds  are connected with dimension $n\geq 2$, and  assumed to be compact while the integrals are taken.\\

By putting \be N^i_j=G^{i\bar{k}}{\dot{\partial}_{\bar{k}}\partial_{j}G}\ee
where $(G^{i\bar k})_{n\times n}=(G_{i\bar j})^{-1}_{n\times n}$ and
\be{\partial}_i:=\frac{\partial}{\partial
z^i},\quad{\partial}_{\bar j}:=\frac{\partial}{\partial \bar
z^j},\quad\dot{\partial}_i:=\frac{\partial}{\partial
v^i},\quad\dot{\partial}_{\bar j}:=\frac{\partial}{\partial \bar
v^j},\ee
the \textit{horizontal vectors} and \textit{vertical covectors} can be defined by
\be \delta_i=\frac{\delta}{\delta z^i}:=\partial_i-N^k_i\dot\partial_k,\ \delta v^i:=dv^i+N^i_kdz^k.\ee
The complexified (co)tangent bundle has the following   horizontal and vertical decomposition
\be T_\mathbb{C}\tilde M=\mathcal{H}\oplus\overline{\mathcal{H}}\oplus\mathcal{V}\oplus\overline{\mathcal{V}},\ \ \ T_\mathbb{C}^*\tilde M=\mathcal{H}^*\oplus\overline{\mathcal{H}^*}\oplus\mathcal{V}^*\oplus\overline{\mathcal{V}^*}\ee
where $\mathcal{H}=\text{span}\{\delta_i\}$, $\mathcal{V}=\text{span}\{\dot\partial_i\}$, $\mathcal{H}^*=\text{span}\{dz^i\}$ and $\mathcal{V}^*=\text{span}\{\delta v^i\}$. Therefore, the operators $\partial$, $\bar\partial$ and $d$ on $\tilde M$ can be decomposed into
\be \partial=\partial_\mathcal{H}+\partial_\mathcal{V},\ \ \overline{\partial}=\overline{\partial}_\mathcal{H}+\overline{\partial}_\mathcal{V},\ \ d=d_\mathcal{H}+d_\mathcal{V}=(\partial_\mathcal{H}+\overline{\partial}_\mathcal{H})+(\partial_\mathcal{V}+\overline{\partial}_\mathcal{V}).\ee
The collection of smooth sections of $(\wedge^p\mathcal{H}^*)\wedge (\wedge^q\overline{\mathcal{H}}^*)\wedge(\wedge^r\mathcal{V}^*)\wedge(\wedge^s\overline{\mathcal{V^*}})$ is denoted by $A^{p,q;r,s}(\tilde M)$, and each element of $A^{p,q;r,s}(\tilde M)$ is called a $(p,q;r,s)$\textit{-form} of $\tilde M$. The elements in $A^{p,q;0,0}(\tilde M)$ are called \textit{horizontal $(p,q)$-forms}. The space of $(l,m)$-forms is clearly $\displaystyle A^{l,m}(\tilde M)=\oplus_{p+r=l,q+s=m}A^{p,q;r,s}(\tilde M)$.

The \textit{K\"ahler form} (fundamental form) of a  complex Finsler metric $G$ is a horizontal $(1,1)$-form defined by
\be \omega_\mathcal{H}=\sqrt{-1}G_{i\bar j}(z,v)dz^i\wedge d\bar z^j.\ee
For a Hermitian metric, $\omega_\mathcal{H}=\sqrt{-1}G_{i\bar j}(z)dz^i\wedge d\bar z^j$ is independent of $v$ and is a $(1,1)$-form living on the base manifold $M$. Generally, $\omega_\mathcal{H}$ lives on $\tilde M$.

\begin{de}[\cite{AP,CS1}] A complex Finsler metric $G$ is said to be  K\"ahler if and only if $d_\mathcal{H}\omega_\mathcal{H}=0$. In this case, $G$ is called a K\"ahler Finsler metric.
\end{de}
The K\"ahler condition is equivalent to the symmetricity of the \textit{Chern-Finsler connection}. In fact, equipping the vertical bundle $\mathcal{V}$  with a inner product $\mathcal{G}$ where $\mathcal{G}(X,Y)=X^i\bar Y^jG_{i\bar j}(z,v)$ for any $X,Y\in \mathcal{V}_{(z,v)}$, the Chern-Finsler connection is just the Hermitian connection of the Hermitian bundle $(\mathcal{V},\mathcal{G})$, and thus the connection 1-forms $(\omega^{i}_{j})$ can be written as
 \be\omega^{i}_{j}=G^{\bar{k}i}\partial G_{j \bar{k}}
 =\Gamma^{i}_{j,k}dz^{k}+C^{i}_{jk}\delta v^{k},\ee
where\be\Gamma^{i}_{j,k}=G^{\bar{l}i}\delta_k G_{j
    \bar{l}},~~~~
        C^{i}_{jk}=G^{\bar{l}i}\dot{\partial}_k G_{j \bar{l}}.
        \ee
The \textit{horizontal torsion} is defined by
\be\theta=\theta^m_{ki}dz^k\wedge dz^i\otimes\delta_m= (\Gamma^m_{i,k}-\Gamma^m_{k,i})dz^k\wedge dz^i\otimes\delta_m.\ee
We call $\vartheta=\vartheta_kdz^k=\theta^m_{km}dz^k$ the \textit{mean horizontal torsion}.

A direct computation gives $$\partial_{\mathcal{H}}\omega_\mathcal{H}=\frac{\sqrt{-1}}{2}(\Gamma^m_{i,k}-\Gamma^m_{k,i})G_{m\bar j}dz^k\wedge dz^i\wedge d\bar z^j.$$

\begin{lem}A complex Finsler metric is    K\"ahler if and only if  $\theta=0$, i.e. $\Gamma^m_{i,k}=\Gamma^m_{k,i}$.
\end{lem}

\section{Integrations on the projectivized bundle}\setcounter{equation}{0}\setcounter{thm}{0}\setcounter{de}{0}\setcounter{prop}{0}\setcounter{coro}{0}\setcounter{lem}{0}
In this section, we will introduce several notions of integration on the {\it projectivized tangent bundle} $\pi: \mathbb{P}(\tilde M)\to M$ where
 $\mathbb{P}(\tilde M):=\tilde M/\mathbb{C}^*$, of which each
fibre is biholomorphic to $\mathbb{CP}^{n-1}$. The complexified bundles $T_\mathbb{C}(\mathbb{P}(\tilde M))$ and $T_\mathbb{C}^*(\mathbb{P}(\tilde M))$ also have  the horizontal and vertical decomposition as (2.4). We shall adopt the same notion $\mathcal{H},\mathcal{V}$ and etc., though the vertical sub-bundle is $(n-1)$-dimensional in this case. The notations $A^{p,q;r,s}(\mathbb{P}(\tilde M))$ and $A^{l,m}(\mathbb{P}(\tilde M))$ have similar definitions with $A^{p,q;r,s}(\tilde M)$ and $A^{l,m}(\tilde M)$ respectively.

Being aware of $G_{i\bar j}(z,\lambda v)=G_{i\bar j}(z, v)$, the K\"ahler form $\omega_\mathcal{H}$ actually lives on $\mathbb{P}(\tilde M)$. We have another $(1,1)$-form $\sqrt{-1}\partial\bar\partial \log G$ which has no mixed part. Considering $v$ as the homogenous coordinate of $\mathbb{P}(\tilde M)$,  it turns out
\be\sqrt{-1}\partial\bar\partial \log G=\omega_\mathcal{V}-\Theta\ee
where
\be \omega_\mathcal{V}=\sqrt{-1}(\log G)_{i\bar j}\delta v^i\wedge\delta\bar
v^{j},\ \ \ (\log G)_{i\bar j}=\dot\partial_i\dot\partial_{\bar j}(\log G)\ee
and $\Theta$ is the \textit{Kobayashi curvature} (\cite{Ko})
\be\Theta=\frac{\sqrt{-1}}{G}K_{i\bar j} dz^i\wedge d\bar z^j,\ \ K_{i\bar j}=-\partial_i\partial_{\bar j}G+G^{k\bar m}(\partial_i\dot\partial_{\bar m}G)(\partial_{\bar j}\dot\partial_kG).\ee
The pull-back $i_z^*\omega_\mathcal{V}=\sqrt{-1}(\log G)_{i\bar j}dv^i\wedge d\bar v^j$ is the Fubini-Study metric on $\mathbb{P}_z:=\pi^{-1}(z)$, where $i_z: \mathbb{P}_z\to \mathbb{P}(\tilde M)$ is the inclusion. Together with $\omega_\mathcal{H}$, the \textit{Sasaki type metric} on $\mathbb{P}(\tilde M)$ is defined as
\be\omega_{\mathbb{P}(\tilde M)}:=\omega_\mathcal{V}+\omega_\mathcal{H}.\ee
The \textit{invariant volume form} can be given by
\be d\mu_{\mathbb{P}(\tilde M)}:=\frac{\omega^{n-1}_\mathcal{V}}{(n-1)!}\wedge\frac{\omega^n_\mathcal{H}}{n!}.\ee

\begin{lem}[\cite{ZZ}] We have $d(\delta_i\lrcorner d\mu_{\mathbb{P}(\tilde M)})=\Gamma^j_{j,i}d\mu_{\mathbb{P}(\tilde M)}$ and its conjugate form $d(\delta_{\bar{i}}\lrcorner d\mu_{\mathbb{P}(\tilde M)})=\overline{\Gamma^j_{j,i}}d\mu_{\mathbb{P}(\tilde M)},$ where ``$\lrcorner$" is the interior derivative.
\end{lem}

Denote $A^{l,m}(M)$  the space of $(l,m)$-forms on $M$. Given $l,m\geq 0$, putting $l^*=l+(n-1), m^*=m+(n-1)$, the \textit{integration along the fibre} is a map $\pi_*: A^{l^*,m^*}(\mathbb{P}(\tilde M))\to A^{l,m}(M)$ which is defined as follows
\be (\pi_*\phi)|_{z}(X_1,\cdots,X_{l},\overline{Y_1},\cdots,\overline{Y_{m}}):=\int_{\mathbb{P}_z}i_z^*\left[\phi(\tilde X_1,\cdots,\tilde X_{l},\overline{\tilde Y_1},\cdots,\overline{\tilde Y_{m}},\cdots)\right]\ee
where $X_i,Y_j\in T_z'M$ and $\tilde X_i,$ $\tilde Y_j$ are their lifts. The RHS of (3.6) is independent of the lifts, and one may use the horizontal ones. Moreover, one can  see that $\pi_*(A^{p,q;r,s}(\mathbb{P}(\tilde M)))=0$ if $r\not=n-1$ or $s\not=n-1$, since $\mathbb{P}_z$ is $(n-1)$-dimensional.

\begin{lem}[cf. \S6 of \cite{BT}] For the bundle $\pi: \mathbb{P}(\tilde M)\to M$, given $\phi\in A(\mathbb{P}(\tilde M))$ and $\alpha\in A(M)$, the integration along the fibre $\pi_*$ satisfies

$(i)$ $d(\pi_*\phi)=\pi_*(d\phi)$;

$(ii)$ $\pi_*((\pi^*\alpha)\wedge\phi)=\alpha\wedge\pi_*\phi$.

\noindent If $M$ is compact in additional, it holds

$(iii)$ $\displaystyle\int_M\alpha\wedge\pi_*\phi=\int_{\mathbb{P}(\tilde M)}\pi^*\alpha\wedge\phi$.
\end{lem}

Applying the above lemma, one can obtain the constancy of the volumes of each fibre which was firstly discovered by R. Yan.

\begin{thm}[\cite{Yan}] Assuming that $(M,G)$ is a complex Finsler manifold, the volume of each fibre $\mathrm{vol}(\mathbb{P}_z):=\pi_*\big(\frac{\omega^{n-1}_\mathcal{V}}{(n-1)!}\big)|_z$ is a constant.
\end{thm}
\noindent\textit{Proof}. Recall $\pi_*\phi=0$ if the vertical part of $\phi$ is not full. Thus
$$\pi_*\omega^{n-1}_\mathcal{V}=\pi_*\left(\sqrt{-1}\partial\bar\partial \log G\right)^{n-1}$$
by (3.1). Hence $$d\left(\mathrm{vol}(\mathbb{P}_z)\right)=d\left(\pi_*\frac{(\sqrt{-1}\partial\bar\partial \log G)^{n-1}}{(n-1)!}\right)=\pi_*\left(d\frac{(\sqrt{-1}\partial\bar\partial \log G)^{n-1}}{(n-1)!}\right)=0.$$
By the connectness of $M$, the volumes of each fibre are constant. \hfill$\Box$

The same technique will give the following rigid result.
\begin{thm}[\cite{Ai}] If $M$ admits a K\"ahler Finsler metric, then it admits a K\"ahler Hermite metric.
\end{thm}
\noindent\textit{Proof}. Let $F$ be a  K\"ahler Finsler metric. Consider the form
$$\omega_M:=\pi_*(\omega_\mathcal{H}\wedge (\sqrt{-1}\partial\bar\partial \log G)^{n-1})=\pi_*(\omega_\mathcal{H}\wedge (\omega_\mathcal{V}-\Theta)^{n-1})=\pi_*(\omega_\mathcal{H}\wedge \omega_\mathcal{V}^{n-1}).$$
Since $d_\mathcal{H}\omega_\mathcal{H}=0$, by Lemma 3.2 We have
\begin{eqnarray*}
  d \omega_M&=& \pi_*\Big((d\omega_\mathcal{H})\wedge (\sqrt{-1}\partial\bar\partial \log G)^{n-1}\Big)\\
  &=& \pi_*\Big((d_\mathcal{V}\omega_\mathcal{H})\wedge (\sqrt{-1}\partial\bar\partial \log G)^{n-1}\Big) \\
  &=& \pi_*\Big((d_\mathcal{V}\omega_\mathcal{H})\wedge \sum_{k=0}^{n-1}C_{n-1}^k\omega_\mathcal{V}^k\wedge(-\Theta)^{n-1-k}\Big).
\end{eqnarray*}
 Recall $\mathrm{dim}_\mathbb{C}\mathbb{P}_z=n-1.$ For $k<n-1$, the vertical part of $(d_\mathcal{V}\omega_\mathcal{H})\wedge \omega_\mathcal{V}^k\wedge(-\Theta)^{n-1-k}$ is not full, thus $\pi_*((d_\mathcal{V}\omega_\mathcal{H})\wedge \omega_\mathcal{V}^k\wedge(-\Theta)^{n-1-k})=0$. For $k=n-1$, the vertical part of $(d_\mathcal{V}\omega_\mathcal{H})\wedge \omega_\mathcal{V}^{n-1}$ overflows. Hence $d\omega_M=0$. One can  deduce the positivity of $\omega_M$ from $\omega_M(X,\overline X)=\pi_*(\omega_\mathcal{H}(\tilde X,\overline {\tilde X})\cdot \omega_\mathcal{V}^{n-1})$.\hfill$\Box$\\

As the end of this section, let us give the definition of the \textit{induced volume form} on $M$.
\begin{de} The induced volume form of $M$ is defined by $d\mu_M:=\pi_*(d\mu_{\mathbb{P}(\tilde M)})$. In other words,
\be\int_M f(z) d\mu_M=\int_{\mathbb{P}(\tilde M)}f(z)d\mu_{\mathbb{P}(\tilde M)}\ee
for any function $f\in C^\infty(M)$.
\end{de}
\noindent \textbf{Remark.} In other literatures, the induced volume form may be divided by a constant and refer to $\frac{1}{\text{vol}(\mathbb{CP}^{n-1})}\pi_*(d\mu_{\mathbb{P}(\tilde M)})$ or $\frac{1}{\text{vol}(\mathbb{P}_z)}\pi_*(d\mu_{\mathbb{P}(\tilde M)})$.

\section{Conformal K\"ahler metrics}\setcounter{equation}{0}\setcounter{thm}{0}\setcounter{de}{0}\setcounter{prop}{0}\setcounter{coro}{0}\setcounter{lem}{0}
Let $G$ be a complex Finsler metric on $M$. A \textit{conformal transformation} of $G$ is a change  $G\mapsto e^fG$ where $f=f(z)$ is a smooth real function on $M$. We denote $e^f G$ by $\hat G$, and the notations of the quantities of $\hat G$ shall wear a hat, e.g. $\hat{\mathcal{H}}$ is the horizontal sub-bundle with respect to $\hat G$ and ${\hat\omega}_{\hat{\mathcal{H}}}$ is the K\"ahler form of $\hat G$. One can easily check
\be \hat G_{i\bar j}=e^fG_{i\bar j},\ \ {\hat\omega}_{\hat{\mathcal{H}}}=e^f{\omega}_{{\mathcal{H}}}\ee
\be \hat N^i_{j}=N^i_{j}+f_jv^i,\ \ \ \hat\Gamma^i_{k,j}=\Gamma^i_{k,j}+f_j\delta^i_k\ee
where $f_j:=\partial_jf.$ Thus
 \be\hat\delta_j=\delta_j-f_jv^i\dot\partial_i,\ \ \ \hat \delta v^i=\delta v^i+v^i\partial f.\ee
Since $v^m\dot\partial_mG_{i\bar j}=0$ by the homogeneity of $G$, we see
\begin{eqnarray}
  \partial_{\hat{\mathcal{H}}}{\hat\omega}_{\hat{\mathcal{H}}} &=& \partial_{\hat{\mathcal{H}}}(\sqrt{-1}e^fG_{i\bar j}dz^i\wedge d\bar z^j) \nonumber\\
   &=&  \hat\delta_k(\sqrt{-1}e^fG_{i\bar j})dz^k\wedge dz^i\wedge d\bar z^j\nonumber\\
   &=&  e^ff_kdz^k\wedge \omega_\mathcal{H}+\sqrt{-1}e^f(\delta_kG_{i\bar j}-f_kv^m\dot\partial_mG_{i\bar j})dz^k\wedge dz^i\wedge d\bar z^j\nonumber\\
   &=&e^f(\partial f\wedge \omega_\mathcal{H}+\partial_\mathcal{H}\omega_\mathcal{H})
\end{eqnarray}
and thus
\be d_{\hat{\mathcal{H}}}{\hat\omega}_{\hat{\mathcal{H}}}=e^f(df\wedge \omega_\mathcal{H}+d_\mathcal{H}\omega_\mathcal{H}).\ee

One can obtain the uniqueness of the K\"ahler Finsler metric in  a conformal class by (4.5). Indeed, a stronger result can be proved.  A Finsler metric is said to be \textit{weakly K\"ahler} if $d_\mathcal{H}\omega_\mathcal{H}(\cdot,\chi,\bar\chi)=0$ where $\chi=v^i\delta_i$ (cf. \cite{AP}). We can show the uniqueness of the weakly K\"ahler Finsler metric in a conformal class.
\begin{thm} In the conformal class of a complex Finsler metric, there exists at most one weakly K\"ahler metric up to homotheties.
\end{thm}

\noindent \textit{Proof}.  By (4.3), one can see that $\hat\chi-\chi$ is vertical. Thus (4.5) gives
$$d_{\hat{\mathcal{H}}}{\hat\omega}_{\hat{\mathcal{H}}}(\cdot,\hat\chi,\bar{\hat\chi})=e^f(df\wedge \omega_\mathcal{H}+d_\mathcal{H}\omega_\mathcal{H})(\cdot,\chi,\bar\chi).$$
If $e^fG$ and $e^gG$ are both weakly K\"ahler, then
$$d(f-g)\wedge\omega_\mathcal{H}(\cdot,\chi,\bar\chi)=0$$
which is equivalent to
  $$(f_i-g_i)G=(f_m-g_m)v^m\dot\partial_iG.$$
Taking the derivative with respect to $\bar v^j$, we get
\begin{eqnarray*}
  (f_{\bar k}-g_{\bar k})\bar v^k(f_m-g_m)v^mG_{i\bar j} &=& (f_{\bar k}-g_{\bar k})\bar v^k(f_i-g_i)\dot\partial_{\bar j}G \\
   &=&  (f_i-g_i) (f_{\bar k}-g_{\bar k})\bar v^k\dot\partial_{\bar j}G \\
   &=&  (f_i-g_i) (f_{\bar j}-g_{\bar j})G.
\end{eqnarray*}
One can easily see that RHS and LHS have different rank unless $d(f-g)=0$.  Therefore,  $e^fG$ and $e^gG$ are homothetic if they are both weakly K\"ahler.\hfill$\Box$\\

At present, let us consider the existence of K\"ahler Finsler metric in the conformal class of a complex Finsler metric. In other words, we shall consider the solvability of the equation \be df\wedge \omega_\mathcal{H}+d_\mathcal{H}\omega_\mathcal{H}=0.\ee

 A Finsler manifold $(M,G)$ is said to be \textit{globally conformal K\"ahler} if and only if there exists a global defined function $f\in C^\infty(M)$ such that $e^fG$ is a K\"ahler Finsler metric. We give the following definition for local solutions.

\begin{de}[cf. \cite{Vai}]A complex Finsler manifold $(M,G)$ is said to be \textit{locally conformal K\"ahler} if and only if there exists an open cover $\{U_\alpha\}$ endowed with smooth functions $f_\alpha:U_\alpha\to \mathbb{R}$ such that $e^{f_\alpha}G$ is a K\"ahler Finsler metric on $U_\alpha$.
\end{de}

By Theorem 4.1, one can see  $d(f_\alpha-f_\beta)=0$ on $U_\alpha\cap U_\beta$ whenever it is nonempty. Thus we obtain a globally defined real 1-form $\varphi\in A^1(M)$ such that $\varphi|_{U_\alpha}=df_\alpha$. Additionally, we have
\be \varphi\wedge \omega_\mathcal{H}+d_\mathcal{H}\omega_\mathcal{H}=0,\ \ \ d\varphi=0.\ee
Such equation was considered by H. Lee \cite{Lee}. Therefore, a real 1-form $\varphi\in A^1(M)$  satisfies (4.7) is called a \textit{Lee form} of $(M,G)$. Thus, if $(M,G)$ is locally conformal K\"ahler, then $(M,G)$ admits a Lee form. Conversely, given a Lee form $\varphi$,  locally  we have $\varphi=df_\alpha$  by Poincar\'{e} Lemma,  and hence $e^{f_\alpha}G$ is a K\"ahler Finsler metric.

\begin{lem} A complex Finsler metric $G$ is locally conformal K\"ahler if and only if $(M,G)$ admits a Lee form.
\end{lem}

On a simply connected manifold, a Lee form is (globally) $d$-exact. Hence, a simply connected, locally conformal K\"ahler manifold is globally conformal K\"ahler.   Moreover, following I. Vaisman \cite{Vai}, we can prove the following rigid theorem.

\begin{thm}Let $(M,G)$ be a compact, locally conformal K\"ahler Finsler manifold. Then $(M,G)$ is globally conformal K\"ahler if and only if $M$ admits a K\"ahler Finsler metric.
\end{thm}
\noindent \textit{Proof}. We prove the sufficiency. Let $\varphi$ be a Lee form of $(M,G)$. We will show that there exist a global function $f\in C^\infty(M)$ such that $\varphi=df.$ Decompose $\varphi$ into $(1,0)$ and $(0,1)$ types
$\varphi=\varphi'+\varphi''$ where $\varphi''=\overline{\varphi'}$. Put $\phi=\sqrt{-1}(\varphi'-\varphi'')$ which is again a real 1-form. We have
\be d\phi=\sqrt{-1}(d\varphi'-d\varphi'')=2\sqrt{-1}d\varphi'=2\sqrt{-1}\bar\partial\varphi'\ee
by $d\varphi=\partial\varphi'+(\bar\partial\varphi'+\partial\varphi'')+\bar\partial\varphi''=0$. Thus $d\phi$ is a real exact $(1,1)$-form.

On the other hand, since $M$ admits a K\"ahler Finsler metric, we have a K\"ahler Hermitian metric  on $M$ by Theorem 3.2. Hence, the $\partial\bar\partial$-lemma holds on the  compact  manifold $M$. Thus, there exists a global real function $f\in C^\infty(M)$ such that
\be\bar\partial\varphi'=\bar\partial\partial f.\ee
Let us consider the metric $\hat G=e^fG.$ Putting $\hat\varphi=\varphi-df$, by (4.5) and (4.7) we have
\be\hat\varphi\wedge{\hat\omega}_{\hat{\mathcal{H}}}+d_{\hat{\mathcal{H}}}{\hat\omega}_{\hat{\mathcal{H}}}=(\varphi-df)\wedge e^f\omega_\mathcal{H}+e^f(df\wedge \omega_\mathcal{H}+d_\mathcal{H}\omega_\mathcal{H})=0.\ee
Therefore $\hat\varphi$ is a Lee form of $(M,\hat G)$.  Write $\hat\varphi=\hat\varphi'+\hat\varphi''$ into $(1,0)$ and $(0,1)$ types. By (4.9) we have
\be \bar\partial\hat\varphi'=\bar\partial(\varphi'-\partial f)=0.\ee
Thus $\hat\varphi'=\hat\varphi_idz^i$ is a holomorphic 1-form. Noting $\overline{\hat\varphi''}=\hat\varphi'$, (4.10) is equivalent to
\be \hat\varphi_i\hat G_{j\bar k}+\hat\delta_i\hat G_{j\bar k}= \hat\varphi_j\hat G_{i\bar k}+\hat\delta_j\hat G_{i\bar k}.\ee
With the help of (2.8), contracting the above equation with $G^{j\bar k}$, one shall reach
\be (n-1)\hat\varphi_i=\hat\Gamma^m_{i,m}-\hat\Gamma^m_{m,i}.\ee
By Lemma 3.1 and (4.11), we finally get
\begin{eqnarray}
  0&=&\int_{\mathbb{P}(\tilde M)}d(\hat\varphi_i\hat G^{i\bar j}\hat\delta_{\bar j}\lrcorner d\hat\mu_{\mathbb{P}(\tilde M)})\nonumber\\
   &=& \int_{\mathbb{P}(\tilde M)}\hat\delta_{\bar j}(\hat\varphi_i\hat G^{i\bar j}) d\hat\mu_{\mathbb{P}(\tilde M)}+\int_{\mathbb{P}(\tilde M)}\hat\varphi_i\hat G^{i\bar j}d(\hat\delta_{\bar j}\lrcorner d\hat\mu_{\mathbb{P}(\tilde M)})\nonumber \\
   &=& \int_{\mathbb{P}(\tilde M)}\hat\varphi_i\hat\delta_{\bar j}(\hat G^{i\bar j}) d\hat\mu_{\mathbb{P}(\tilde M)}+\int_{\mathbb{P}(\tilde M)}\hat\varphi_i\hat G^{i\bar j}{\overline{\hat\Gamma^{m}_{m,j}}} d\hat\mu_{\mathbb{P}(\tilde M)}\nonumber \\
   &=& \int_{\mathbb{P}(\tilde M)}\hat\varphi_i\hat G^{i\bar j}(-\overline{\hat\Gamma^{m}_{j,m}}+\overline{\hat\Gamma^{m}_{m,j}})d\hat\mu_{\mathbb{P}(\tilde M)}\nonumber\\
   &=&-(n-1)\int_{\mathbb{P}(\tilde M)}\hat\varphi_i\hat G^{i\bar j}\overline{\hat\varphi_j}d\hat\mu_{\mathbb{P}(\tilde M)}
\end{eqnarray}
which implies $\hat\varphi'=\hat\varphi_idz^i=0$. Hence, $\hat\varphi=0$ and  $\varphi=df.$\hfill$\Box$\\

Theorem 4.2 tells us that the equation (4.6) is globally solvable if and only if it is locally solvable, if the compact manifold $M$ admits a K\"ahler Finsler metric. Recalling the definitions of the horizontal torsion $\theta$, the equation (4.6) can be expressed in the form
\be f_k\delta^i_j-f_j\delta^i_k=\Gamma^i_{k,j}-\Gamma^i_{j,k}=\theta^i_{jk}.\ee
The trace of (4.15) gives
\be (n-1)f_j=-\vartheta_j\ee
where $\vartheta_j$'s are the components of the mean horizontal torsion $\vartheta$.

\begin{thm} Let $M$ be a compact manifold admitting a K\"ahler Finsler metric. Then, a complex Finsler metric $G$ on $M$ is globally conformal K\"ahler if and only if

\noindent $(i)$ the horizontal torsion is reducible $\theta^i_{jk}=\frac{1}{n-1}(\vartheta_j\delta^i_k-\vartheta_k\delta^i_j)$
 where $\vartheta_j=\theta^m_{jm}$;

\noindent $(ii)$ and $d(\vartheta+\bar{\vartheta})=0.$
\end{thm}
\noindent \textit{Proof}. One can easily get the necessity by (4.15) and (4.16). Conversely,
$$d(\vartheta+\overline{\vartheta})=\partial\vartheta+\partial\bar{\vartheta}+\bar\partial\vartheta+\bar\partial\bar{\vartheta}=0$$
implies ${\bar\partial}_\mathcal{V}\bar\vartheta=(\dot\partial_{\bar j}\bar\vartheta_i)\delta \bar v^j\wedge d\bar z^i=0$. Thus $\vartheta_i=\vartheta_i(z)$ is independent of $v$, and $\vartheta$ must be  a 1-form living on the base manifold $M$. Then by the Poincar\'{e} Lemma, $(n-1)df=-(\vartheta+\bar\vartheta)$ is locally solvable on $M$, which implies $(n-1)\partial f=-\vartheta$ locally. Together with $(i)$, we  get (4.15). Finally, (4.15) is globally solvable by Theorem 4.2.\hfill$\Box$

\section{Total holomorphic curvature}\setcounter{equation}{0}\setcounter{thm}{0}\setcounter{de}{0}\setcounter{prop}{0}\setcounter{coro}{0}\setcounter{lem}{0}

In this section, we will consider the total holomorphic curvature in the conformal classes. Let us recall the definition of the curvature forms.
The curvature forms $\Omega^i_j:=\bar\partial\omega^i_j$ of the Chern-Finsler connection can be divided into four parts, namely, $h\bar h$-, $v\bar h$-, $h\bar v$ and $v\bar v$-curvatures. By (2.7), the $h\bar h$-curvature has the form
\be h\bar h\text{-component of }\Omega^i_j=R^i_{j,k\bar m}dz^k\wedge d\bar z^m=(-\delta_{\bar m}\Gamma^i_{j,k}-C^i_{js}\delta_{\bar m}N^s_k)dz^k\wedge d\bar z^m.\ee
Putting $R_{j\bar l,k\bar m}:=G_{i\bar l}R^i_{j,k\bar m}$, a direct computation gives (cf.\cite{AP})
\be R_{j\bar l,k\bar m}v^j\bar v^l=K_{k\bar m}\ee
where $K_{k\bar m}$ is the Kobayashi curvature given in (3.3). The \textit{holomorphic curvature} is defined by
\be K(z,v):=\frac{1}{G^2} R_{j\bar l,k\bar m}v^j\bar v^lv^k\bar v^m=\frac{1}{G^2}K_{i\bar j}v^i\bar v^j.\ee
We define
the {\it total holomorphic curvature}  of $(M,G)$ by setting
\be {\mathcal K}(G)=\int_{\mathbb{P}(\tilde M)}K(z,v)d\mu_{\mathbb{P}(\tilde M)}.\ee
In order to consider the above functional in the \textit{volume preserved conformal class} \be [G]=\big\{e^fG\,|\, f\in C^\infty(M), \mathrm{vol}(M,e^fG)=\mathrm{vol}(M,G)\big\},\ee
 let us give a divergence lemma.

 \begin{lem}Given $\alpha=\alpha_idz^i\in A^{1,0;0,0}(\mathbb{P}(\tilde M))$, we have
 \be d(\alpha_iG^{i\bar j}\delta_{\bar j}\lrcorner d\mu_{\mathbb{P}(\tilde M)})=G^{i\bar j}(\alpha_{i,\bar j}+\alpha_i\bar\vartheta_j)d\mu_{\mathbb{P}(\tilde M)}\ee
 \be d\Big(\frac{1}{G}\alpha_iv^i\cdot\bar \chi\lrcorner d\mu_{\mathbb{P}(\tilde M)}\Big)=\frac{1}{G}\Big(\alpha_{i,\bar j}v^i\bar v^j+\alpha_iv^i\cdot\bar\vartheta_j\bar v^j\Big)d\mu_{\mathbb{P}(\tilde M)}\ee
 and their conjugate forms, where $\alpha_{i,\bar j}:=\delta_{\bar j}\alpha_i$ and $\chi=v^j\delta_j$.
 \end{lem}
\noindent \textit{Proof}. The proof of (5.6) is similar to (4.14). For (5.7), applying $\delta_{\bar j}G=0$, $v^j\Gamma^i_{j,k}=N^i_k$ and Lemma 3.1, we get
\begin{eqnarray*}
  d\Big(\frac{1}{G}\alpha_iv^i\cdot\bar\chi\lrcorner d\mu_{\mathbb{P}(\tilde M)}\Big) &=& \delta_{\bar j}\Big(\frac{1}{G}\alpha_i v^i\bar v^j\Big) d\mu_{\mathbb{P}(\tilde M)}+\frac{1}{G}\alpha_iv^i\bar v^jd(\delta_{\bar j}\lrcorner d\mu_{\mathbb{P}(\tilde M)})\\
   &=& \frac{1}{G}\Big(\alpha_{i,\bar j}v^i\bar v^j-\alpha_iv^i\cdot\overline{N^j_j}+ \alpha_iv^i \cdot\overline{v^j\Gamma^m_{m,j}}\Big)d\mu_{\mathbb{P}(\tilde M)}\\
  &=&  \frac{1}{G}\Big(\alpha_{i,\bar j}v^i\bar v^j+\alpha_iv^i\cdot\bar\vartheta_j\bar v^j\Big)d\mu_{\mathbb{P}(\tilde M)}.
\end{eqnarray*}
The conjugate forms of (5.6) and (5.7) are obviously true.\hfill$\Box$

At present, let us give the relations of the curvatures of two conformal related metrics. Putting $\hat G=e^{f(z)}G$, by (4.1)-(4.3), we get $\hat\omega_{\hat{\mathcal{H}}}=e^f\omega_{\mathcal{H}}$ and
\be \hat\omega_{\hat{\mathcal{V}}}=\sqrt{-1}(\log G)_{i\bar j}(\delta v^i+v^i\partial f)\wedge (\delta \bar v^j+\bar v^j\bar\partial f)=\omega_{\mathcal{V}}\ee
where we use $(\log G)_{i\bar j}v^i=(\log G)_{i\bar j}\bar v^j=0$ for the last equality. Thus the fibre volume $\text{vol}(\mathbb{P}_z)$ is a conformal invariant, and
\be d{\hat\mu}_{\mathbb{P}(\tilde M)}=e^{nf}d\mu_{\mathbb{P}(\tilde M)},\ \ \  d{\hat\mu}_{M}=e^{nf}d\mu_{M}.\ee
Recalling (3.3), one can obtain
\be \hat K_{i\bar j}=e^f(K_{i\bar j}-f_{i\bar j}G)\ee
where $f_{i\bar j}=\partial_i\partial_{\bar j}f.$ Invariantly, it says
\be \hat\Theta=\Theta-\sqrt{-1}\partial\bar\partial f.\ee
Moreover, one can get
\be \hat K=e^{-f}\left(K-\frac{1}{G}f_{i\bar j}v^i\bar v^j\right).\ee

Now by considering a family of conformal deformations $e^{f(t,z)}G$ with the initial date $f(0,z)=0$, one can find

\begin{eqnarray}
  \frac{d}{dt}\mathcal{K}(e^{f}G) &=&\frac{d}{dt} \int_{\mathbb{P}(\tilde M)}e^{(n-1)f}\left(K-\frac{1}{G}f_{i\bar j}v^i\bar v^j\right)d\mu_{\mathbb{P}(\tilde M)} \nonumber\\
  &=&  \int_{\mathbb{P}(\tilde M)}(n-1)f' e^{(n-1)f}\left(K-\frac{1}{G}f_{i\bar j}v^i\bar v^j\right)d\mu_{\mathbb{P}(\tilde M)}\nonumber\\
  &&\hskip1cm-\int_{\mathbb{P}(\tilde M)}\frac{1}{G}e^{(n-1)f}f'_{i\bar j}v^i\bar v^jd\mu_{\mathbb{P}(\tilde M)},
\end{eqnarray}
where $f'=\frac{\partial}{\partial t}f$.

Denoting $f'(0,z):=\nu(z)$, and substituting $f(0,z)=0$ into (5.13), it turns out
\be \left.\frac{d}{dt}\mathcal{K}(e^{f}G)\right|_{t=0}=\int_{\mathbb{P}(\tilde M)}\Big((n-1)\nu K-\frac{1}{G}\nu_{i\bar j}v^i\bar v^j\Big)d\mu_{\mathbb{P}(\tilde M)}\ee
Taking $\alpha_i=\nu_i$ in (5.7), we get
\be \int_{\mathbb{P}(\tilde M)}\frac{1}{G}\nu_{i\bar j}v^i\bar v^jd\mu_{\mathbb{P}(\tilde M)}=-\int_{\mathbb{P}(\tilde M)}\frac{1}{G}\nu_{i}v^i\bar\vartheta_j\bar v^jd\mu_{\mathbb{P}(\tilde M)}.\ee
Then taking $\alpha_i=\nu\vartheta_i$ in (5.7), its conjugate form gives
\be \int_{\mathbb{P}(\tilde M)}\frac{1}{G}\nu_{i}v^i\bar\vartheta_j\bar v^jd\mu_{\mathbb{P}(\tilde M)}=- \int_{\mathbb{P}(\tilde M)}\frac{1}{G}\nu(|\vartheta_iv^i|^2+\overline{\vartheta_{i,\bar j}v^i\bar v^j})d\mu_{\mathbb{P}(\tilde M)}.\ee
Note that $\nu_{i\bar j}v^i\bar v^j$ is real, we obtain

\be \left.\frac{d}{dt}\mathcal{K}(e^{f}G)\right|_{t=0}=\int_{\mathbb{P}(\tilde M)}\nu\Big((n-1)K-\frac{1}{G}(|\vartheta_iv^i|^2+\mathfrak{Re}(\vartheta_{i,\bar j}v^i\bar v^j))\Big)d\mu_{\mathbb{P}(\tilde M)}.\ee
At this point, let us define the \textit{mean holomorphic curvature} $\kappa$ by
\be \kappa={\pi_*\big(Kd\mu_{\mathbb{P}(\tilde M)}\big)}\big/{\pi_*\big(d\mu_{\mathbb{P}(\tilde M)}\big)}\ee
which is a real function on $M$, and call
\be \kappa_{\vartheta}={\pi_*\Big(\Big(K-\frac{1}{G(n-1)}\big(|\vartheta_iv^i|^2+\mathfrak{Re}(\vartheta_{i,\bar j}v^i\bar v^j)\big)\Big)d\mu_{\mathbb{P}(\tilde M)}\Big)}\big/{\pi_*\big(d\mu_{\mathbb{P}(\tilde M)}\big)}\ee
the $\vartheta$-\textit{mean holomorphic curvature}.
By (5.7), one can see
$$\int_{\mathbb{P}(\tilde M)}\Big(|\vartheta_iv^i|^2+\mathfrak{Re}(\vartheta_{i,\bar j}v^i\bar v^j)\Big)d\mu_{\mathbb{P}(\tilde M)}=0.$$
Recalling $\pi_*(d\mu_{\mathbb{P}(\tilde M)})=d\mu_M$, we obtain various representations of $\mathcal{K}(G)$
\be \int_M\kappa_\vartheta d\mu_M=\int_M\kappa d\mu_M=\int_M\pi_*(Kd\mu_{\mathbb{P}(\tilde M)})=\int_{\mathbb{P}(\tilde M)}Kd\mu_{\mathbb{P}(\tilde M)}=\mathcal{K}(G).\ee
Since $\pi_*(\nu\phi)=\nu\pi_*(\phi)$ for any $\phi\in A(\mathbb{P}(\tilde M))$, the   formula (5.17) becomes
\be \left.\frac{d}{dt}\mathcal{K}(e^{f}G)\right|_{t=0}=(n-1)\int_M\nu\kappa_\vartheta d\mu_M.\ee
Assuming $e^{f}G$ in the volume preserved class $[G]$, we have
$$0=\frac{d}{dt}\int_Me^{nf}d\mu_M=n\int_Mf'e^{nf}d\mu_M.$$
At $t=0$, it reads as
$$0=\int_M\nu d\mu_M.$$
Thus a critical point shall satisfies
\be \int_M\nu\kappa_\vartheta d\mu_M=0\ \ \ \text{where}\ \ \int_M\nu d\mu_M=0.\ee
Denoting the average  $\bar{\kappa}_\vartheta=\frac{1}{\mathrm{vol}(M)}\int_M\kappa_\vartheta d\mu_M,$ it is equivalent to
\be \int_M\nu(\kappa_\vartheta-\bar\kappa_\vartheta) d\mu_M=0\ \ \ \text{where}\ \ \int_M\nu d\mu_M=0.\ee
Taking $\nu=\kappa_\vartheta-\bar\kappa_\vartheta$, it becomes
\be \int_M(\kappa_\vartheta-\bar\kappa_\vartheta)^2 d\mu_M=0.\ee

\begin{thm} A metric $G$ is a critical point of $\int_M\kappa_\vartheta d\mu_M$ in its volume preserved conformal class $[G]$ if and only if $\kappa_\vartheta=const$. If $G$ is a K\"ahler Finsler metric, then $\kappa=const.$
\end{thm}
Particularly, a K\"ahler Finsler metric with constant holomorphic curvature is critical in the volume preserved conformal class. Next, let us consider the \textit{stability of a critical K\"ahler Finsler metric}. The second variation is

\begin{eqnarray}
\frac{d^2}{dt^2}\mathcal{K}(e^fG)&=&\int_{\mathbb{P}(\tilde M)}(n-1)f'' e^{(n-1)f}\left(K-\frac{1}{G}f_{i\bar j}v^i\bar v^j\right)d\mu_{\mathbb{P}(\tilde M)}\nonumber \\
&&+\int_{\mathbb{P}(\tilde M)}(n-1)^2f'f' e^{(n-1)f}\left(K-\frac{1}{G}f_{i\bar j}v^i\bar v^j\right)d\mu_{\mathbb{P}(\tilde M)}\nonumber\\
  &&-\int_{\mathbb{P}(\tilde M)}\frac{2(n-1)}{G}f' e^{(n-1)f}f'_{i\bar j}v^i\bar v^jd\mu_{\mathbb{P}(\tilde M)}\nonumber\\
  && -\int_{\mathbb{P}(\tilde M)}\frac{1}{G}e^{(n-1)f}f''_{i\bar j}v^i\bar v^jd\mu_{\mathbb{P}(\tilde M)}
\end{eqnarray}
where $f'=\frac{\partial}{\partial t}f$ and $f''=\frac{\partial^2}{\partial t^2}f$. At $t=0$, denoting $f''(0,z)=\psi(z)$ and recalling $f(0,z)=0$ and $f'(0,z)=\nu(z)$, we get

\begin{eqnarray}
\left.\frac{d^2}{dt^2}\mathcal{K}(e^fG)\right|_{t=0}&=&\int_{\mathbb{P}(\tilde M)}(n-1)\psi Kd\mu_{\mathbb{P}(\tilde M)}+\int_{\mathbb{P}(\tilde M)}(n-1)^2\nu^2Kd\mu_{\mathbb{P}(\tilde M)}\nonumber\\
  &&-\int_{\mathbb{P}(\tilde M)}\frac{2(n-1)}{G}\nu \nu_{i\bar j}v^i\bar v^jd\mu_{\mathbb{P}(\tilde M)} -\int_{\mathbb{P}(\tilde M)}\frac{1}{G}\psi_{i\bar j}v^i\bar v^jd\mu_{\mathbb{P}(\tilde M)}.\nonumber\\
  &&
\end{eqnarray}
Since $G$ is a K\"ahler Finsler metric, the torsion $\vartheta$ vanishes. Taking $\alpha_i=\nu\nu_i$ in (5.7), we get
\be \int_{\mathbb{P}(\tilde M)}\frac{1}{G}\nu \nu_{i\bar j}v^i\bar v^jd\mu_{\mathbb{P}(\tilde M)}=-\int_{\mathbb{P}(\tilde M)}\frac{1}{G}\nu_{\bar j} \nu_{i}v^i\bar v^jd\mu_{\mathbb{P}(\tilde M)}\ee
while taking $\alpha_i=\psi_i$, it leads to
\be \int_{\mathbb{P}(\tilde M)}\frac{1}{G}\psi_{i\bar j}v^i\bar v^jd\mu_{\mathbb{P}(\tilde M)}=0.\ee
By defining  a \textit{induced Hermitian metric $h$}
\be h^{i\bar j}:={\pi_*\Big(\frac{2}{G}v^i\bar v^jd\mu_{\mathbb{P}(\tilde M)}\Big)}\big/{\pi_*\big(d\mu_{\mathbb{P}(\tilde M)}\big)},\ee
the equation (5.26) becomes
\be \left.\frac{d^2}{dt^2}\mathcal{K}(e^fG)\right|_{t=0}=(n-1)\int_M\left(h^{i\bar j}\nu_i\nu_{\bar j}+(\psi+(n-1)\nu^2)\kappa\right)d\mu_M.\ee
Let us recall
\be 0=\left.\frac{d^2}{dt^2}\right|_{t=0}\int_Me^{nf}d\mu_M=\int_Mn(\psi+n\nu^2)d\mu_M.\ee
Thus, by the constancy of $\kappa$, finally we have

\be\left.\frac{d^2}{dt^2}\mathcal{K}(e^fG)\right|_{t=0}=(n-1)\int_M\left(h^{i\bar j}\nu_i\nu_{\bar j}-\nu^2\kappa\right)d\mu_M\ee
where $\int_M\nu d\mu_M=0$. We call $G$ a \textit{stable} critical metric of $\mathcal{K}$   if the above second variation is nonnegative.

\begin{thm} In a volume preserved conformal class, a critical K\"ahler Finsler metric of the functional $\int_M\kappa d\mu_M$ is stable if and only if
the constant mean holomorphic curvature satisfies $\kappa\leq\lambda_1^h$, where $\lambda_1^h$ is the first eigenvalue of the Hermitian Laplacian of the metric measure space $(M,h,d\mu_M)$ given by
\be \lambda_1^h:=\inf\left\{\frac{\int_Mh^{i\bar j}\phi_i\phi_{\bar j}d\mu_M}{\int_M \phi^2d\mu_M}\,\left|\,\phi\in C^\infty(M),\ \int_M\phi d\mu_M=0\right.\right\}.\ee
\end{thm}

\section{Total holomorphic Ricci curvature}\setcounter{equation}{0}\setcounter{thm}{0}\setcounter{de}{0}\setcounter{prop}{0}\setcounter{coro}{0}\setcounter{lem}{0}
In this section, we will consider the Ricci curvature of a complex Finsler metric.
The \textit{holomorphic Ricci curvature} of $G$ is defined as
\be Ric(z,v)=\frac{1}{G}G^{i\bar j}R_{k\bar m,i\bar j}v^k\bar v^m=\frac{1}{G}G^{i\bar j}K_{i\bar j}.\ee
Kobayashi introduced an analogous quantity for complex Finsler vector bundles in \cite{Ko2}, and named it the mean curvature.

The \textit{total holomorphic Ricci curvature} of $(M,G)$ is given by
\be {\mathcal R}(G)=\int_{\mathbb{P}(\tilde M)}Ric(z,v)d\mu_{\mathbb{P}(\tilde M)}.\ee
Denoting $\hat G=e^fG$ again, one can deduce
\be\hat Ric=e^{-f}(Ric-G^{i\bar j}f_{i\bar j})\ee
from (5.10).
By a similar calculation of \S5, we have
\be \left.\frac{d}{dt}\mathcal{R}(e^{f}G)\right|_{t=0}=\int_{\mathbb{P}(\tilde M)}(n-1)\nu Ric\,d\mu-\int_{\mathbb{P}(\tilde M)} G^{i\bar j}\nu_{i\bar j}d\mu_{\mathbb{P}(\tilde M)}.\ee
Taking $\alpha=\nu_idz^i,$  one can deduce from (5.6) that
\be- \int_{\mathbb{P}(\tilde M)} G^{i\bar j}\nu_{i\bar j} d\mu_{\mathbb{P}(\tilde M)}=\int_{\mathbb{P}(\tilde M)} G^{i\bar j}\nu_{i}\bar\vartheta_j d\mu_{\mathbb{P}(\tilde M)}.\ee
Taking $\alpha=\nu\vartheta_idz^i$, the conjugate form of (5.6) gives
\be\int_{\mathbb{P}(\tilde M)} G^{i\bar j}\nu_{i}\bar\vartheta_j d\mu_{\mathbb{P}(\tilde M)}=-\int_{\mathbb{P}(\tilde M)}\nu G^{i\bar j}(\vartheta_i\bar\vartheta_j+\overline{\vartheta_{j,\bar i}}) d\mu_{\mathbb{P}(\tilde M)}.\ee
Since the expression is real, we obtain

\be\left.\frac{d}{dt}\mathcal{R}(e^{f}G)\right|_{t=0}=\int_{\mathbb{P}(\tilde M)}\nu\Big((n-1)Ric-\big(\|\vartheta\|_G^2+\mathfrak{Re}(\vartheta_{i,\bar j}G^{i\bar j})\big)\Big)d\mu_{\mathbb{P}(\tilde M)}.\ee
Let us define the \textit{mean holomorphic Ricci curvature} $\rho$ by
\be \rho={\pi_*\big(Ric\,d\mu_{\mathbb{P}(\tilde M)}\big)}\big/{\pi_*\big(d\mu_{\mathbb{P}(\tilde M)}\big)}\ee
which is again a real function on $M$. We call
\be \rho_{\vartheta}={\pi_*\Big(\big(Ric-\frac{1}{(n-1)}\big(\|\vartheta\|_G^2+\mathfrak{Re}(\vartheta_{i,\bar j}G^{i\bar j})\big)\big)d\mu_{\mathbb{P}(\tilde M)}\Big)}\big/{\pi_*\big(d\mu_{\mathbb{P}(\tilde M)}\big)}\ee
the $\vartheta$-\textit{mean holomorphic Ricci curvature}.
By (5.6), one can see
\be \int_{\mathbb{P}(\tilde M)}\Big(\|\vartheta\|_G^2+\mathfrak{Re}(\vartheta_{i,\bar j}G^{i\bar j})\Big)d\mu_{\mathbb{P}(\tilde M)}=0\ee
and thus
\be \int_M\rho_\vartheta d\mu_M=\int_M\rho\, d\mu_M=\int_M\pi_*(Ric\,d\mu_{\mathbb{P}(\tilde M)})=\int_{\mathbb{P}(\tilde M)}Ric\,d\mu_{\mathbb{P}(\tilde M)}=\mathcal{R}(G).\ee
By the definition of $\rho_\vartheta$, the first variation formula (6.7) becomes
\be \left.\frac{d}{dt}\mathcal{R}(e^{f}G)\right|_{t=0}=(n-1)\int_M\nu\rho_\vartheta d\mu_M.\ee

\begin{thm} A metric $G$ is a critical point of $\int_M\rho_\vartheta d\mu_M$ in its volume preserved conformal class $[G]$ if and only if $\rho_\vartheta=const$. If $G$ is a K\"ahler Finsler metric, then $\rho=const.$
\end{thm}

Let $G$ be a critical K\"ahler Finsler metric. We shall give its second variation formula. Similarly to \S5, we have
\begin{eqnarray}
\left.\frac{d^2}{dt^2}\mathcal{R}(e^fG)\right|_{t=0}&=&\int_{\mathbb{P}(\tilde M)}(n-1)\psi Ric\,d\mu_{\mathbb{P}(\tilde M)}+\int_{\mathbb{P}(\tilde M)}(n-1)^2\nu^2Ric\,d\mu_{\mathbb{P}(\tilde M)}\nonumber\\
  &&-\int_{\mathbb{P}(\tilde M)}2(n-1)\nu \nu_{i\bar j}G^{i\bar j}d\mu_{\mathbb{P}(\tilde M)} -\int_{\mathbb{P}(\tilde M)}\psi_{i\bar j}G^{i\bar j}d\mu_{\mathbb{P}(\tilde M)}\nonumber\\
  &=&\int_{M}(n-1)\psi \rho\, d\mu_{M}+\int_{M}(n-1)^2\nu^2\rho\, d\mu_{M}\nonumber\\
  &&+\int_{\mathbb{P}(\tilde M)}2(n-1)\nu_{\bar j} \nu_{i}G^{i\bar j}d\mu_{\mathbb{P}(\tilde M)}.
\end{eqnarray}
Let us define another \textit{induced Hermitian metric $g$ by}
\be g^{i\bar j}:={\pi_*\Big(2G^{i\bar j}d\mu_{\mathbb{P}(\tilde M)}\Big)}\big/{\pi_*\big(d\mu_{\mathbb{P}(\tilde M)}\big)}.\ee
By (5.31), we have
\be \left.\frac{d^2}{dt^2}\mathcal{R}(e^fG)\right|_{t=0}=(n-1)\int_M(g^{i\bar j}\nu_i\nu_{\bar j}-\nu^2\rho)d\mu_M.\ee
where $\int_M\nu d\mu_M=0$.
Finally, we can  state the \textit{stability} of a critical K\"ahler Finsler metric of the functional $\mathcal{R}$.

\begin{thm} In a volume preserved conformal class, a critical K\"ahler Finsler metric of the functional $\mathcal{R}=\int_M\rho d\mu_M$ is stable if and only if
the constant mean holomorphic Ricci curvature satisfies $\rho\leq\lambda_1^g$, where $\lambda_1^g$ is the first eigenvalue of the Hermitian Laplacian  of the metric measure space $(M,g,d\mu_M)$ defined by
\be \lambda_1^g:=\inf\left\{\frac{\int_Mg^{i\bar j}\phi_i\phi_{\bar j}d\mu_M}{\int_M \phi^2d\mu_M}\,\left|\,\phi\in C^\infty(M),\ \int_M\phi d\mu_M=0\right.\right\}.\ee
\end{thm}
We adopt Kobayashi's notion of Finsler Einstein bundles (\cite{Ko2}) and give the following definition of K\"ahler Finsler metrics.
\begin{de} A K\"ahler Finsler metric with constant holomorphic Ricci curvature is called a K\"ahler  Einstein Finsler  metric.
\end{de}
By this definition, one can immediately get the following corollary.
\begin{coro} A K\"ahler Einstein Finsler metric with non-positive holomorphic Ricci curvature  is a stable critical point of $\mathcal{R}$ in its volume preserved conformal class.
\end{coro}

\section{A Yamabe type problem}\setcounter{equation}{0}\setcounter{thm}{0}\setcounter{de}{0}\setcounter{prop}{0}\setcounter{coro}{0}\setcounter{lem}{0}

In this section, we shall study the existence of complex Finsler metrics with constant $\rho_\vartheta$ in the volume preserved conformal class $[G]$.  Through the variational approach (cf. \cite{CZ,LP}), we can get the existence of metrics with $\rho_\vartheta=const.$

Customary, write the conformal change in the form $\hat G=\phi^{\frac{2}{n-1}}G$, where $\phi$ is a positive function and $n$ is the complex dimension of $M$. Consider the following Yamabe type functional
\be \mathfrak{R}(\phi)=\frac{1}{\mathrm{vol}^{1-\frac{1}{n}}(M,\phi^{\frac{2}{n-1}}G)}\mathcal{R}(\phi^{\frac{2}{n-1}}G).\ee
Using Lemma 5.1, (6.9) and (6.14), we have
\begin{eqnarray}
  \nonumber &&\int_{\mathbb{P}(\tilde M)}\hat Ric\,d\hat\mu_{\mathbb{P}(\tilde M)}\\
  \nonumber &=&\int_{\mathbb{P}(\tilde M)}(\phi^2Ric+\frac{2}{n-1}(\phi_i\phi_{\bar j}-\phi_{i\bar j}\phi)G^{i\bar j})d\mu_{\mathbb{P}(\tilde M)} \\
  \nonumber &=& \int_{\mathbb{P}(\tilde M)}(\phi^2Ric+\frac{1}{n-1}(4G^{i\bar j}\phi_i\phi_{\bar j}-\phi^2\|\vartheta\|_G^2-\phi^2\mathfrak{Re}(\vartheta_{i,j}G^{i\bar j}))d\mu_{\mathbb{P}(\tilde M)}  \\
   &=&\int_M(\frac{2}{n-1}g^{i\bar j}\phi_i\phi_{\bar j}+\phi^2\rho_\vartheta)d\mu_M.
\end{eqnarray}
In the real expression, $g^{i\bar j}\phi_i\phi_{\bar j}$ is $\frac{1}{4}\|d\phi\|_g^2$,  thus the Yamabe type functional (7.1) is of the form
\be  \mathfrak{R}(\phi)=\frac{\int_M(\frac{1}{2(n-1)}\|d\phi\|_g^2+\phi^2\rho_\vartheta)d\mu_M}{\left(\int_M\phi^{\frac{2n}{n-1}}d\mu_M\right)^{\frac{n-1}{n}}}.\ee
By the H\"{o}lder's inequality, one can get $\mathfrak{R}(\phi)\geq-\left(\int_M|\rho_\vartheta|^nd\mu_M\right)^{1/n}$, thus we can defined a \textit{conformal invariant} as
\be Y(G)=\inf_{0<\phi\in C^\infty(M)}\mathfrak{R}(\phi).\ee
The \textit{energy} of $\phi$ is given by
\be E(\phi)=\int_M(\frac{1}{2(n-1)}\|d\phi\|_g^2+\phi^2\rho_\vartheta)d\mu_M\ee
and the $L^q$-norm is defined as $\|\phi\|_q=\left(\int_M|\phi|^qd\mu_M\right)^{1/q}.$
By putting $p=\frac{2n}{n-1}$, we have
\be \mathfrak{R}(\phi)=\frac{E(\phi)}{\|\phi\|_p^2}.\ee
Since $C^\infty(M)$ is dense in the Sobolev space $W^{1,2}(M)$, $\mathfrak{R}(|\phi|)\leq\mathfrak{R}(\phi)$ and $\mathfrak{R}(\lambda\phi)=\mathfrak{R}(\phi)$ for $\lambda>0$,  we see
$$Y(G)=\inf_{\phi\in W^{1,2}}\mathfrak{R}(\phi)=\inf_{\|\phi\|_p=1}E(\phi).$$
The Euler-Lagrangian equation of the minimizer with $\|\phi\|_p=1$ is
\be L\phi:=\frac{1}{2(n-1)}\Delta_g\phi+\frac{1}{2(n-1)}\langle d\phi,d\log\tau\rangle_g-\phi\rho_\vartheta=-Y(G)\phi^{p-1}.\ee
where $\Delta_g$ is the Laplacian of the induced Hermitian metric $g$ and $\tau=\frac{d\mu_M}{d\mu_g}$.

Note that the real dimension of $M$ is $m=2n$, therefore $p=\frac{2n}{n-1}=\frac{2m}{m-2}$ is the critical exponent of the Sobolev embedding theorem. Following Yamabe, let us consider the disturbed functional
\be \mathfrak{R}_t(\phi)=\frac{E(\phi)}{\|\phi\|_t^2},\ \ \ \  2\leq t\leq p=\frac{2n}{n-1},\ee
whose infimum is denoted by $Y_t$. The Euler-Lagrangian equation of the minimizer of $\mathfrak{R}_t(\phi)$ with $\|\phi\|_t=1$ is
\be L\phi=-Y_t\phi^{t-1}.\ee
By the regularity theory, for any $t<p$ there exists a smooth and positive minimizer $\phi_t$ of $\mathfrak{R}_t$ with $\|\phi_t\|_t=1$ (cf. Lemma 5.2 in \cite{CZ} or Proposition 4.2 in \cite{LP}). In other  words, for any $2\leq t<p$ we have a smooth and positive function $\phi_t$ satisfies
\be L\phi_t=-Y_t\phi^{t-1}_t.\ee

At this point, we shall consider the limit when $t\to p^+.$ Henceforth, let us assume the initial metric $G$ has unit volume $\mathrm{vol}(M, G)=1$.
\begin{lem}[cf. Lemma 4.3 in \cite{LP}]Given $\mathrm{vol}(M,G)=1$,we have

$(1)$ if $Y_p<0$, then $\limsup_{t\to p^-}Y_t\leq Y_p=Y(G)$;

$(2)$ if $Y_p\geq 0$, then $\lim_{t\to p^-}Y_t=Y_p=Y(G)$.

\end{lem}
As we did in \cite{CZ}, let us introduce another \textit{conformal invariant}
\be C(G)=\sup_{x\in M}\left[\frac{d\mu_g}{d\mu_M}\right]^{\frac{1}{n}}.\ee
 By Definition 3.1 and (6.14), when $G$ is Hermitian, it holds
$C(G)=\frac{1}{2\mathrm{vol}(\mathbb{CP}^{n-1})^{1/n}}$ which can be considered as the normalizing factor of $Y(G)$. Then we have a Sobolev inequality.

\begin{lem}[cf.  Lemma 5.4 in \cite{CZ}] Let $(M,G)$ be a compact complex Finsler manifold. Then for any $\epsilon>0$, there exists $C_\epsilon$ such that
\be \|w\|_p^2\leq\frac{(1+\epsilon)C(G)}{\sigma_{2n}}\int_M\|dw\|^2_{g}d\mu_M+C_\epsilon\int_Mw^2d\mu_M\ee
where $\sigma_{2n}$ is the best Sobolev constant on $\mathbb{R}^{2n}$ satisfies
\be \sigma_{2n}\left(\int_{\mathbb{R}^{2n}}|f|^pdx\right)^{\frac{2}{p}}\leq\int_{\mathbb{R}^{2n}}\|df\|^2dx.\ee
\end{lem}

\noindent \textit{Proof}. Recalling $\tau=d\mu_{M}/d\mu_g$, let us put $\tilde g_{i\bar j}=\tau^{\frac{1}{n}}g_{i\bar j}$. It turns out $d\mu_{\tilde g}=d\mu_M$ and thus (cf. Theorem 2.3 in \cite{LP})
$$\|w\|_p^2\leq\frac{(1+\epsilon)}{\sigma_{2n}}\int_M\|dw\|^2_{\tilde g}d\mu_{\tilde g}+C_\epsilon\int_Mw^2d\mu_{\tilde g}.$$
We can deduce (7.12) from
$\|dw\|^2_{\tilde g}=\tau^{-1/n}\|dw\|_g^2\leq C(G)\|dw\|_g^2.$ \hfill$\Box$

According to Lemma 7.1-7.2, by a similar argument of Proposition 4.4 in \cite{LP}, one can obtain the following uniform $L^{p_0}$ estimate.
\begin{lem} If $Y(G)\cdot C(G)<\frac{\sigma_{2n}}{2n-2}$, then there exists $t_0<p$ and $p_0>p$ such that $\phi_t(t_0\leq t<p)$ are uniformly bounded in $L^{p_0}$.
\end{lem}

Finally, the regularity theory gives $\{\phi_t\}$ are uniformly  bounded in $C^{2,\alpha}(M)$.  Then $\phi_{t_i}\to \phi$ in $C^2(M)$ for some $t_i\to p$, and the limit gives
$-L\phi\leq Y(G)\phi^{p-1}$, $\|\phi\|_p=1$  and $\mathfrak{R}(\phi)\leq Y(G)$. Hence $\mathfrak{R}(\phi)=Y(G)$ by the definition of $Y(G)$. Moreover, the minimizer $\phi$ satisfies $-L\phi=Y(G)\phi^{p-1}$ and then $\phi$ is smooth and positive.

\begin{thm}If $Y(G)\cdot C(G)<\frac{\sigma_{2n}}{2n-2}$, then there exists a smooth positive function $\phi$ such that $\mathfrak{R}(\phi)=Y(G).$ In this case, there exists a metric $\hat G$ in the conformal class $[G]$ such that $\hat\rho_\vartheta=const.$
\end{thm}

As the end, we shall give the following upper bound theorem.

\begin{thm} For any compact complex Finsler  manifold $(M,G)$, it holds $Y(G)\cdot C(G)\leq \frac{2\sigma_{2n}}{n-1}$.
\end{thm}

\noindent\textit{Proof}.  The proof is similar to the real case we given in \cite{CZ}. Recall that $m=2n$ is the real dimension of $M$. It is well-known that the function
\be u_\epsilon:=\left(\frac{\epsilon}{\epsilon^2+r^2}\right)^{\frac{m-2}{2}},\ \ \ r=|x|,\ \ \epsilon>0\ee
achieve the best Sobolev constant on the Euclidean space $\mathbb{R}^m$ and satisfies
$$\partial_ru_\epsilon=-(m-2)\frac{r}{\epsilon^2+r^2}u_\epsilon,\ \Delta_{\mathbb{R}^m}u_\epsilon=-m(m-2)u_\epsilon^{p-1}$$
which imply
\begin{align}
 \int_{B(R)-B(\rho)}|du_\epsilon|^2dx=&m(m-2)\int_{B(R)-B(\rho)} u_\epsilon^pdx\nonumber\\
 &+(2-m)\omega_{m-1}\epsilon^{m-2}\left[\frac{R^m}{(\epsilon^2+R^2)^{m-1}}-\frac{\rho^m}{(\epsilon^2+\rho^2)^{m-1}}\right]
\end{align}
 where $B(R)=\{x:|x|<R\}$ and $\omega_{m-1}=\mathrm{vol}(\mathbb{S}^{m-1})$.
Hence the Sobolev constant satisfies
\be\sigma_{2n}=\sigma_{m}=\frac{\int_{\mathbb{R}^m}|du_\epsilon|^2dx}{\left(\int_{\mathbb{R}^m}u_\epsilon^pdx\right)^{\frac{2}{p}}}=m(m-2)\left(\int_{\mathbb{R}^m}u_\epsilon^pdx\right)^{\frac{2}{m}}.\ee
Moreover, we have
\be\int_{B(\rho)}|du_\epsilon|^2dx<m(m-2)\int_{B(\rho)} u_\epsilon^pdx< \sigma_{m}\left(\int_{B(\rho)} u_\epsilon^p dx\right)^{\frac{2}{p}}\ee
and
\be  \int_{B(\rho)} u_\epsilon^p dx =\omega_{m-1}\int_0^\rho  \left( \frac{\epsilon}{\epsilon^2+r^2}\right)^m r^{m-1}dr
= \omega_{m-1}\int_0^{\rho/\epsilon} \frac{t^{m-1}}{(1+t^2)^m}dt. \ee

Let $\eta=\eta(r)$ be a radial cutoff function on $\mathbb{R}^m$, such that $ 0\leq \eta\leq1$, $\eta|_{B(1)}=1$, $\eta|_{\mathbb{R}^m-B(2)}=0$, and $|d\eta|=|\partial_r\eta|\leq 2.$ Putting $\eta_\rho:=\eta(\frac{r}{\rho})$ for $\rho>0$, we have $ 0\leq \eta_\rho\leq1$, $\eta|_{B(\rho)}=1$, $\eta|_{\mathbb{R}^m-B(2\rho)}=0$, and $|d\eta_\rho|=|\partial_r\eta_\rho|\leq \frac{2}{\rho}$. Consider the test function
$\varphi:=\eta_\rho u_\epsilon$ for $\epsilon<<\rho$.

Recall $\tau=d\mu_{M}/d\mu_g$ and $\tilde g=\tau^{\frac{1}{n}}g=[d\mu_{M}/d\mu_g]^{\frac{1}{n}}g$. Let us pick a point $x_0\in M$ such that $C(G)=\sup_{x\in M}\tau^{-1/n}(x)=\tau^{-1/n}(x_0),$ and take a normal coordinate  system of $\tilde g$ centered at $x_0$. By the continuity,  we have
 $$\tau^{1/n}(x)\leq \frac{1}{C(G)}+\delta(\rho),\ \ x\in B(2\rho)$$
 where $\delta(\rho)\to0$ when $\rho\to0$.
 Suppose $2\rho$ is less than the injectivity radius of $x_0$ with respect to $\tilde g$. The test function
$\varphi=\eta_\rho u_\epsilon$ can be considered as a globally defined function on $M$. We will give the estimate of $\mathfrak{R}(\varphi)=\frac{E(\varphi)}{\|\varphi\|_p^2}$.

Applying the relations between $\tilde g$ and $g$, we have
$$ E(\varphi)\leq  \frac{1}{2n-2}\int_M \tau^{\frac{1}{n}}\|d\varphi\|^2_{\tilde g}d\mu_{\tilde g}+c_1\int_M\varphi^2 d\mu_{\tilde g}.$$
  Assume $(1-c_2|x|)dx\leq d\mu_{\tilde g}\leq (1+c_2|x|)dx$  in $B(2\rho)$.  By the H\"{o}lder inequality and (7.18), one gets the estimate
$$\int_M\varphi^2 d\mu_{\tilde g} \leq (1+2c_2\rho)\int_{B(2\rho)}u_\epsilon^2dx\leq c_3 \left(\int_{B(2\rho)}  u_\epsilon^pdx\right)^{\frac{2}{p}}\rho^2\leq c_4\rho^2.$$
Next, we give the estimate of the term
$$\int_{B(2\rho)} \tau^{1/n}\|d\varphi\|^2_{\tilde g}d\mu_{\tilde g}\leq \left(\frac{1}{C(G)}+\delta(\rho)\right)\int_{B(2\rho)}  \|d\varphi\|^2_{\tilde g}d\mu_{\tilde g}.$$
Since the space is  locally Euclidean, one can obtain
\begin{align}
\int_M\|d\varphi\|^2_{\tilde g}d\mu_{\tilde g}&\leq(1+2c_2\rho) \int_{B(2\rho)}|\partial_r\varphi|^2dx\nonumber\\
&=(1+2c_2\rho) \left[\int_{B(\rho)}|\partial_ru_\epsilon|^2dx+\int_{B(2\rho)-B(\rho)}|\partial_r(\eta_\rho u_\epsilon)|^2dx\right].\nonumber
\end{align}
The first term can be estimated by (7.17). For the second term, we see from (7.15) that
\begin{align*}
\int_{B(2\rho)-B(\rho)}|\partial_r(\eta_\rho u_\epsilon)|^2dx&\leq\frac{8}{\rho^2} \int_{B(2\rho)-B(\rho)}  u_\epsilon^2dx+2\int_{B(2\rho)-B(\rho)}  |\partial_ru_\epsilon|^2dx\\
&\leq c_5 \left(\int_{B(2\rho)-B(\rho)}  u_\epsilon^pdx\right)^{\frac{2}{p}}+c_5\int_{B(2\rho)-B(\rho)}   u_\epsilon^pdx\\
&\ \ \ \ +c_5\rho^{2-m}\epsilon^{m-2}.
\end{align*}
Being aware of (7.18), we see that
$$(1+2c_2\rho)\int_{B(2\rho)-B(\rho)}|\partial_r(\eta_\rho u_\epsilon)|^2dx\leq\frac{c_6\epsilon^{m-2}}{\rho^{m-2}}.$$
On the other hand, for any $\epsilon<\rho<\frac{1}{2c_2}$, it holds
\be\left(\int_{M} \varphi^pd\mu_{M}\right)^{\frac{2}{p}}=\left(\int_{M} \varphi^pd\mu_{\tilde g}\right)^{\frac{2}{p}}\geq  (1-c_2\rho)^{\frac{2}{p}}\left(\int_{B(\rho)}  u_\epsilon^pdx\right)^{\frac{2}{p}}\geq c_{7}.\ee
Together with (7.16)-(7.19), we reach
 $$\mathfrak{R}(\varphi)\leq \left(\frac{1}{C(G)}+\delta(\rho)\right)\left[\frac{(1+2c_2\rho)}{(1-c_2\rho)^{\frac{2}{p}}}\frac{\sigma_{2n}}{2n-2}+\frac{c_6\epsilon^{n-2}}{c_{7}\rho^{n-2}}\right]+\frac{c_1c_4}{c_7}\rho^2.$$
By letting $\epsilon\to 0$ and $\rho\to0$, we see $Y(G)\leq  \frac{1}{C(G)}\cdot \frac{\sigma_{2n}}{2n-2}.$\hfill$\Box$\\

\noindent\textbf{Remark}. The same procedure can be used to study the existence of metrics with constant $\kappa_\vartheta$.

\newpage
\small{

\noindent Bin Chen

School of Mathematical Sciences, Tongji University

Shanghai, 200092, P. R. China

chenbin@tongji.edu.cn\\

\noindent Yibing Shen

School of Mathematical Sciences, Zhejiang University

Hangzhou, 310027, P. R. China

yibingshen@zju.edu.cn\\

\noindent Lili Zhao

School of Mathematical Sciences, Shanghai Jiao Tong University

Shanghai, 200240, P. R. China

zhaolili@sjtu.edu.cn
}

\end{document}